\documentclass{amsart}
\usepackage{amssymb,amscd,verbatim}
\usepackage{enumerate}

\usepackage{graphicx}

\newcommand{\NN}{\mathbb N}

\newcommand{\RR}{\mathbb R}

\newcommand{\ZZ}{\mathbb Z}
\newcommand{\CI}{{\mathcal C}^{\infty}}
\newcommand{\CIc}{{\mathcal C}^{\infty}_{\text{c}}}
\newcommand\pa{{\partial}}
\newcommand{\loc}{\operatorname{loc}}

\newcommand{\maC}{\mathcal C}

\newcommand{\maH}{\mathcal H}

\newcommand{\ie}{{\em i.e., }}
\newtheorem{theorem}{Theorem}[section]
\newtheorem{proposition}[theorem]{Proposition}

\theoremstyle{definition}
\newtheorem{definition}[theorem]{Definition}
\theoremstyle{remark}

\numberwithin{figure}{section}
\numberwithin{table}{section}

\newcommand\Kond[2]{{\mathcal K}^{#1}_{#2}}

\newcommand\OO[1]{\Omega^{(#1)}}

\newcommand{\cyl}{\maC}
\newcommand{\Bacuta}{{B\u{a}cu\c{t}\u{a}}}

\author[C. \Bacuta]{Constantin \Bacuta}
\email{bacuta@math.udel.edu}
\author[V. Nistor]{Victor Nistor}
\email{nistor@math.psu.edu}
\author[L. Zikatanov]{Ludmil T. Zikatanov}
\email{ltz@math.psu.edu}
\address{Pennsylvania State University,
Math. Dept., University Park, PA 16802}
\thanks{The authors were supported in part by the NSF grant DMS
02-09497.  Victor Nistor was also partially supported by NSF grant
DMS 02-00808.}

\begin{document}
\date{\today}

\title[Boundary value problems]{Regularity and well posedness for
the Laplace operator on polyhedral domains}

\begin{abstract}
We announce a well-posedness result for the Laplace equation in
weighted Sobolev spaces on polyhedral domains in $\RR^n$ with
Dirichlet boundary conditions. The weight is the distance to the
set of singular boundary points. We give a detailed sketch of the
proof in three dimensions.
\end{abstract}

\maketitle \tableofcontents

\section*{Introduction}

Let $\Omega \subset \RR^n$ be an open set. Consider the boundary
value problem
\begin{equation}\label{eq.BVP}
    \Delta u = f, \quad  u\vert_{\pa \Omega} = g,
\end{equation}
where $\Delta$ is the Laplace operator. For $\Omega$ smooth and
bounded, this boundary value problem has a unique solution $u \in
H^{s+2}(\Omega)$ depending continuously on $f \in H^s(\Omega)$ and
$g \in H^{s+3/2}(\pa \Omega)$, $s \ge -1$. See the books of Evans
\cite{Evans} or Taylor \cite{Taylor1} for a proof of this basic
and well known result.

It is also well known that this result does not extend to
non-smooth domains $\Omega$. A deep analysis of the difficulties
that arise for $\pa \Omega$ Lipschitz is contained in the papers
of Babu\v{s}ka and Guo \cite{BabuGuo}, \Bacuta, Bramble, and Xu
\cite{bacutaBX}, Jerison and Kenig \cite{JerisonKenig}, Kenig
\cite{Kenig}, Mitrea and Taylor \cite{MitreaTaylor}, Verchota
\cite{Verchota}, and others (see the references in the
aforementioned papers). Results specific to polyhedral domains are
contained in the papers of Costabel \cite{Costabel}, Dauge
\cite{Dauge, Dauge1}, Elschner \cite{Elschner1}, Kondratiev
\cite{Kondratiev67}, Mazya and Rossman \cite{MazyaRossmann} and
others. Good references are also the monographs \cite{Grisvard1,
KMRossmann, MazyaBoth}.

In this paper, we consider the boundary value problem
\eqref{eq.BVP} when $\Omega$ is a \emph{bounded  polyhedron}
in $\RR^n$, and Poisson's equation $\Delta u = f$ is
replaced by a strongly elliptic system. Let us denote by
$\OO{n-2}$ the set of points $p \in \pa \Omega$ such $\pa \Omega$
is not smooth in a neighborhood of $p$ and by $\eta_{n-2}(x)$ the
distance from a point $x \in \Omega$ to the set $\OO{n-2} \subset
\pa \Omega$ of non-smooth boundary points of $\Omega$.

We shall work in the {\em weighted Sobolev spaces}
\begin{equation}\label{eq.def.wSsp0}
    \Kond{\mu}{a}(\Omega) = \{u \in L^2_{\loc}(\Omega),\,
    \eta_{n-2}^{|\alpha|-a} \pa^\alpha u \in L^2(\Omega),
    \text{ for all } |\alpha| \le \mu\}, \quad \mu \in \ZZ_+,
\end{equation}
which we endow with the induced Hilbert space norm. We shall call
these spaces the {\em Babu\v{s}ka--Kondratiev} spaces. A similar
definition (\ie using the same weight, see Definition
\ref{eq.def.wSsp}) yields the Babu\v{s}ka--Kondratiev (or weighted
Sobolev) spaces $\Kond{\mu}{a}(\pa \Omega)$, $\mu \in \ZZ_+$. The
spaces $\Kond{s}{a}(\pa \Omega)$, $s \in \RR_+$ are defined by
interpolation. The Babu\v{s}ka--Kondratiev spaces are closely related
to weighted Sobolev spaces on non-compact manifolds. See the works of
Erkip and Schrohe \cite{ErkipSchrohe} and Grubb \cite{Grubb}, for
related results on boundary value problems on non-compact manifolds
and, more generally, on the analysis on non-compact manifolds. Here is
our main result for the Laplace equation on a polyhedral domain.

\begin{theorem}\label{theorem.Main0}\
Let $\Omega \subset \RR^n$ be a bounded, polyhedral domain and
$\mu \in \ZZ_+$. Then there exists $\eta > 0$ such that the
boundary value problem \eqref{eq.BVP} has a unique solution $u \in
\Kond{\mu+1}{a+1}(\Omega)$ for any $f \in
\Kond{\mu-1}{a-1}(\Omega)$, any $g \in \Kond{\mu+1/2}{a+1/2}(\pa
\Omega)$, and any $|a|< \eta$. This solution depends continuously
on $f$ and $g$. If $\mu = 0$, this solutions is the solution of
the associated variational problem.
\end{theorem}

The proof can be carried out, without much change, to yield the same
result for strongly elliptic, strictly positive systems on curvilinear
polyhedral domains. The complications are mostly of topological
nature, so we shall discuss this result in \cite{BNZ2}. The analytic
part of the proof is however the same both for polyhedral domains and
for curvilinear polyhedral domains, therefore the reader interested
mostly in analysis will benefit from the simplified account included
in this paper.

We now describe the contents of the sections of the paper in more
detail. The first section introduces the weighted Sobolev spaces (also
called the {\em Babu\v{s}ka--Kondratiev} spaces) that appear in our
main result, Theorem \ref{theorem.Main0}. The second section contains
a statement of three intermediate results:\ a Hardy--Poincar\'{e}
inequality, a regularity result, and a trace theorem. A proof of proof
of the Hardy--Poincar\'e type inequality in dimensions $n = 3$ is
given in the third section. A sketch of the proof of our main result
is given in Section \ref{sec.Main}. This proof is based on the three
intermediate results mentioned above. The last two intermediate
results are particular cases of some results proved in \cite{AIN},
provided that we show that polyhedral domains fit into the framework
of Lie domains developed in that paper. This is however highly
non-trivial in higher dimensions and leads to topological and
geometric complications that will be treated in detail in \cite{BNZ2}
in the more general framework of curvilinear polyhedral domains.

We thank Bernd Ammann, Ivo Babu\v{s}ka, Wolfgang Dahmen, Alexandru
Ionescu, and Daniel T\u{a}taru for useful discussions. We also
thank Johnny Guzman for pointing the reference \cite{Kondratiev70}
to us. The second named author would like to thank Institute Henri
Poincar\'e in Paris for its hospitality while parts of this work
were being completed.

\section{Sobolev spaces\label{sec.SP}}

We introduce in this section the Babu\v{s}ka--Kondratiev (or weighted
Sobolev) spaces $\Kond{\mu}{a}(\Omega)$, $\Kond{s}{a}(\pa \Omega)$,
$\mu \in \ZZ$ and $s \in \RR$, for the case when $\Omega$ is a {\em
straight polyhedron} (straight polyhedra are defined below by
induction).  Throughout this paper $\Omega$ will be an open
set. Recall that we denoted by $\eta_{n-2}(x)$ the distance from a
point $x \in \RR^n$ to $\OO{n-2} \subset \pa \Omega$, the set of
singular boundary points of $\Omega$.

Below, by an {\em affine space} we shall denote the translation of
a subspace of a vector space $V$. We define the concept of a
straight polyhedron of dimension $n$ by induction. A subset
$\Omega$ of an affine line is a {\em straight polyhedron of
dimension $1$} if it is a finite union of open intervals (bounded
or not) on a line such that $\pa \Omega = \pa \overline{\Omega}$.
An open subset $\Omega \subset V$ with finitely many connected
components of an affine space $V$ of dimension $n \ge 2$ will be
called a {\em straight polyhedron of dimension $n$} if $\pa \Omega
= \pa \overline{\Omega}$ and there exist disjoint straight
polyhedra $D_j \subset \pa \Omega$ of dimension $n-1$ such that
$\pa \Omega = \cup \overline{D}_j$.

The condition $\pa \Omega = \pa \overline{\Omega}$ is equivalent
to the fact that $\Omega$ is the interior of its closure
$\overline{\Omega}$. This condition is designed to rule out the
case when $\Omega$ lies of {\em both} sides of its boundary. To
deal with this case, as well as with more general domains, we need
the concept of a ``curvilinear polyhedral domain,'' which will be
discussed in \cite{BNZ2}. A simple example of a polyhedron is the
interior of the convex hull of a finite set of points in $\RR^n$,
provided that this set is not empty.

Let $\Omega$ be a straight polyhedron in an affine space $V$.  For
simplicity we shall take $V=\RR^n$. Let $f$ be a continuous function
on $\Omega$, $f > 0$ on the interior of $\Omega$. We define the
$\mu$th Sobolev space with weight $f$ (and index $a$) by
\begin{equation}\label{eq.def.wSsp}
    \Kond{\mu}{a, f}(\Omega) = \{u \in L^2_{\loc}(\Omega),\,
    f^{|\alpha|-\,a}\pa^\alpha u \in L^2(\Omega),
    \text{ for all } |\alpha| \le \mu \}\,,\quad \mu \in \ZZ_+.
\end{equation}
The norm on $\Kond{\mu}{a, f}(\Omega)$ is $\|u\|_{\Kond{\mu}{a,
f}(\Omega)}^2 := \sum_{|\alpha| \le \mu} \|f^{|\alpha|-\,a}
\pa^\alpha u\|_{L^2(\Omega)}^2$.


\begin{definition}\label{def.main.S}\
We let $\Kond{\mu}{a}(\Omega) = \Kond{\mu}{a, f}(\Omega)$ and
$\Kond{\mu}{a}(\pa \Omega) = \Kond{\mu}{a, f}(\pa \Omega)$, where
$f = \eta_{n-2}$ is the distance to $\OO{n-2}$.
\end{definition}

For example, $\Kond{0}{0}(\Omega) = L^2(\Omega)$. For $\Omega$ a
polygon in the plane, $\eta_{n-2}(x) = \eta_{0}(x)$ is the
distance from $x$ to the vertices of $\Omega$ and the resulting
spaces $\Kond{\mu}{a}(\Omega)$ are the spaces introduced by
Kondratiev \cite{Kondratiev67}. Let us notice that we define both
Sobolev spaces $\Kond{\mu}{a}(\Omega)$ and $\Kond{\mu}{a}(\pa
\Omega)$ using {\em the same} weight function.

If $\mu \in \NN = \ZZ_+ \smallsetminus \{0\}$, we define
$\Kond{-\mu}{-a}(\Omega)$ to be the dual of
\begin{equation}
    \stackrel{\circ}{\Kond{\mu}{a}}(\Omega) := \Kond{\mu}{a}(\Omega)
    \cap \{\pa_\nu^j u \vert_{\pa \Omega} = 0, \ j = 0, 1, \ldots,
    \mu-1\}
\end{equation}
with pivot $\Kond{0}{0}(\Omega)$. Since $\CIc(\Omega))$ is dense
in $\stackrel{\circ}{\Kond{\mu}{a}}(\Omega)$ by Theorem 3.4 of
\cite{AIN}, an equivalent definition of the space
$\Kond{-\mu}{-a}(\Omega)$, $-\mu \in \NN$, is as follows. First
define for any $u \in \CI(\Omega)$
\begin{equation}
    \|u\|_{\Kond{-l}{-a}(\Omega)} = \sup
    \frac{|(u, v)|}{\|v\|_{\Kond{\mu}{a}(\Omega)}}\,,
    \quad 0 \neq v \in \CIc(\Omega).
\end{equation}
Then we let $\Kond{-l}{-a}(\Omega)$ to be the completion of the
space of smooth functions $u$ on $\Omega$ that are such that
$\|u\|_{\Kond{-l}{-a}(\Omega)} < \infty$. The spaces
$\Kond{s}{a}(\pa \Omega)$, with $s \not \in \ZZ$, are defined by
complex interpolation.

The following result can be proved in small dimensions directly
using spherical or polar coordinates. In higher dimensions, it
follows using also the result of \cite{AIN}, and it will be dealt
with in \cite{BNZ2}.

\begin{proposition}\label{prop.map}\
Let $P$ be a differential operator of order $m$ on $\Omega$ with
smooth coefficients. Then $P$ maps $\Kond{\mu}{a}(\Omega)$ to
$\Kond{\mu-m}{a-m}(\Omega)$ continuously, for any admissible
weight $h$ and any $\mu \in \ZZ$. Moreover, the resulting family
$h^{-\lambda} P h^{\lambda} : \Kond{\mu}{a}(\Omega) \to
\Kond{\mu-m}{a-m}(\Omega)$ of bounded operators depends
continuously on $\lambda$.
\end{proposition}

\section{Three intermediate results\label{sec.PM}}

We now state in the three main intermediate results needed for the
proof of our main result, Theorem \ref{theorem.Main0}, namely, a
Hardy--Poincar\'e type inequality (Theorem \ref{thm.Poincare}, a
regularity result for polyhedra (Theorem \ref{thm.reg}), and a theorem
on the general properties of the trace map between weighted Sobolev
spaces (Theorem \ref{thm.Sobolev}).

\begin{theorem}\label{thm.Poincare}\
There exists a constant $\kappa_{\Omega} > 0$, depending only on
$\Omega$, such that
\begin{equation*}
    \|u\|_{\Kond{1}{1}(\Omega)}^2
    \le \kappa_{\Omega} \int_{\Omega} |\nabla u(x)|^2\, dx\,, \quad
    dx = dx_1 dx_2 \ldots dx_n,
\end{equation*}
for any function $u \in H^1_{\loc}(\Omega)$ such that $u
\vert_{\pa \Omega} = 0$.
\end{theorem}

The regularity result, stated next, is of independent interest.

\begin{theorem}\label{thm.reg}\
Assume that $\Delta u \in \Kond{\mu-1}{-1}(\Omega)$ and $u
\vert_{\pa \Omega} \in \Kond{\mu+1/2}{1/2}(\pa \Omega)$, $\mu \in
\ZZ_+$, for some $u \in \Kond{1}{1}(\Omega)$. Then $u \in
\Kond{\mu+1}{1}(\Omega)$ and
\begin{equation}
    \|u\|_{\Kond{\mu+1}{1}(\Omega)} \le
    C\Big(\|\Delta u\|_{\Kond{\mu-1}{-1}(\Omega)}
    + \|u\|_{\Kond{0}{1}(\Omega)} +
    \|u\vert_{\pa \Omega}\|_{\Kond{\mu+1/2}{1/2}(\pa \Omega)}\Big).
\end{equation}
\end{theorem}

We shall need also the following result on weighted Sobolev spaces,
which generalizes the well known results on Sobolev spaces on
domains with smooth boundary. Let $\CIc(\Omega)$ be the space of
compactly supported functions on the open set $\Omega$.

\begin{theorem}\label{thm.Sobolev}\
The restriction $\CIc(\overline{\Omega} \smallsetminus \OO{n-2})
\ni u \to u\vert_{\pa \Omega} \in \CIc(\pa \Omega \smallsetminus
\OO{n-2})$ extends to a continuous, surjective map
\begin{equation*}
    \Kond{\mu}{a}(\Omega) \to \Kond{\mu-1/2}{a-1/2}(\pa \Omega),
    \quad \mu \ge 1.
\end{equation*}
Moreover, $\CIc(\Omega)$ is dense in the kernel of this map if
$\mu = 1$.
\end{theorem}

Theorems \ref{thm.reg} and \ref{thm.Sobolev} will follow from
Theorems 3.4 and 3.7 of \cite{AIN}, once we will have identified
our weighted Sobolev spaces on $\Omega$ with the Sobolev spaces
introduced in \cite{AIN}. This will be done, in the more general
setting of curvilinear polyhedral domains in \cite{BNZ2}. The
proofs of the quoted results from \cite{AIN} is to reduce to the
case of a half-space using a suitable partition of unity. The
construction of this partition of unity is, in turn, based on the
geometric framework of Lie manifolds introduced in \cite{aln1}.

Let us give only a brief hint of the role of Lie algebras of
vector fields (and Lie manifolds) in the study of weighted Sobolev
spaces on a polyhedron. There exits an explicit smooth function
$r_{\Omega}$ on $\Omega$ such that $r_{\Omega}$ is equivalent to
$\eta_{n-2}$ (\ie $r_\Omega/\eta_{n-2}$ is bounded from above and
bounded away from $0$) and, moreover,
$r_{\Omega}^t\Kond{\mu}{a}(\Omega) = \Kond{\mu}{a+t}(\Omega)$.
This function is constructed as follows. Let $g_0$ be the
Euclidean metric. We shall define the metrics $g_{k+1}$ and the
functions $\tilde\rho_{k}$, $k \ge 1$, as follows. Let $\rho_{k}$
be the distance to the faces of dimension $k$ of
$\overline{\Omega}$ in the metric $g_k$ and let $\tilde\rho_{k}$
be a smooth function that coincides with $\rho_k$ when $\rho_k$ is
small and otherwise satisfies $\rho_k/2 \le \tilde\rho_k \le
\rho_k$. We then let $g_{k+1} = \tilde\rho_k^{-2} g_k$. Finally,
we define $\eta_{n-2} = \tilde \rho_0 \tilde \rho_1 \ldots \tilde
\rho_{n-2}$. The vector fields that we are interested are of the
form $f(x) r_{\Omega}X$, where $X$ is a vector field on a
neighborhood of $\Omega$ and $f$ is a function that is smooth in
suitable generalized spherical coordinates. The set of such vector
fields is closed under the Lie bracket of vector fields. See
\cite{BNZ2} for more details.

\section{A Poincar\'e type inequality\label{sec.PI}}

The rest of this section is devoted to a proof of the
Hardy--Poincar\'e type inequality stated in Theorem
\ref{thm.Poincare} in dimension is $n=3$. A proof by induction for
arbitrary $n$ is included in \cite{BNZ2}. That proof requires,
however, the use of curvilinear polyhedral domains on the unit
sphere, which explains why it is convenient to use domains more
general than the straight polyheral ones in higher dimensions.

\subsection{Proof of the Hardy--Poincar\'{e} type inequality for
$n=3$} The idea of the proof is to cover the domain $\Omega$ with
open sets $\cyl$ on which the integration simplifies and we can
use the usual Poincar\'{e} inequality. Then we add the
corresponding inequalities.

We shall write $dV = dx dy dz$ for the volume element. Note that
$\eta_{n-2} = \eta_1$, since we have fixed our dimension.

Let us consider the apparently weaker inequality
\begin{equation}\label{eq.newpoincare}
    \|u\|_{\Kond{0}{1}(\cyl)}^2 :=
    \int_{\cyl} \frac{|u(x)|^2}{\eta_{1}(x)^2}\; dx
    \le C \int_{\cyl} |\nabla u(x)|^2\, dx, \quad u = 0
    \text{ on } \cyl \cap \pa \Omega.
\end{equation}
For $\cyl = \Omega$, as in the case of smooth bounded domains,
this inequality is immediately seen to be equivalent to our
result. Hence we shall concentrate on proving this inequality for
suitable $\cyl$, including $\cyl = \Omega$. In fact, the proof of
our inequality for $\cyl = \Omega$ will be obtained by adding
certain analogues of the inequality \eqref{eq.newpoincare} for
suitable domains $\cyl$.

Assume that $u$ is a smooth function on $\overline{\Omega}$ with
$u\vert_{\pa \Omega} = 0$. We  further consider two small enough
positive numbers $\epsilon > \delta > 0$, depending only on
$\Omega$, such that the following three properties are satisfied:

{\bf 1.}\ For any edge $e$ of $\Omega$, we consider the right
cylindrical domain $Cil_{e}$ of radius $\delta$ whose axis of
symmetry is the line containing the edge $e$ and whose bases
intersect $e$ at distance $\epsilon$ from its two vertices. (These
two basis are orthogonal to $e$, by assumption.) We assume that
$\epsilon$ and $\delta$ were chosen small enough so that the
domain $\Omega_e = \Omega \cap Cil_{e}$ can be characterized, in
suitable cylindrical coordinates, by
\begin{equation*}
    \Omega_e= \{ (r, \theta, z) |\; 0 < r < r_{e}, \ 0 < \theta <
    \theta_{e},\ 0 < z < z_{e} := |e|-2\epsilon \},
\end{equation*}
where $|e|$ is the length of the edge $e$, and $\eta_1 = r$ on
$\Omega_e$. In these cylindrical coordinates, the edge $e$ is on
the $z$-axis (in particular, $r = 0$ on $e$).

{\bf 2.}\ For any vertex $v$ and any edge $e$ containing $v$, we
consider the right conical domain $Con_{v,e}$ with vertex $v$ and
base the same with one of the bases of $Cil_{e}$ (the one which
closer to the vertex $v$) and whose symmetry axis is the line
containing the edge $e$. We assume that $\epsilon$ and $\delta$
were chosen small enough so that domain $\Omega_{v,e} = \Omega
\cap Con_{v,e}$ can be characterized in cylindrical coordinates by
\begin{equation*}
    \Omega_{v,e}= \{ (r, \theta, z) \vert\; 0 < r < z \frac
    {\delta}{\sqrt{\epsilon^2 + \delta^2}},\ 0 < \theta <
    \theta_e,\ 0 < z < \epsilon \}
\end{equation*}
and $\eta_1 = r$ on $\Omega_{v, e}$.
%
%
{\bf 3.}\ Let $B(v, t)$ be the open ball of radius $t$ centered at
$v$. For any vertex $v$ of $\Omega$, the domain $\Omega_v = \Omega
\cap B(v, 2\epsilon)$ can be characterized in (suitable) spherical
coordinates centered at $v$ by
\begin{equation*}
    \Omega_v = \{ (\rho, \omega) \vert\, 0 < \rho < 2\epsilon, \
    \omega \in \omega_v\},
\end{equation*}
where $B(v,r)$ is the three dimensional ball centered   at $v$ and
of radius $r$,  $\omega_v$ is a "polygonal region" on the unit
sphere $S^2 \subset R^2$, and $\rho=0$ corresponds to $v$.

We shall now prove \eqref{eq.newpoincare} for $\cyl$ one of the
domains $\Omega_e$ or $\Omega_{v, e}$. We need to stress at this
poin the crucial importance of the relation $\eta_1 = r$ on these
domains.

Let $W_a$ be the angle $0 < r < a$ and $0 < \alpha$. Let us next
prove first the inequality \eqref{eq.newpoincare} when $\cyl =
\Omega_{e}$, that is, when $\cyl$ is the cylindrical domain
described in cylindrical coordinates $(r, \theta, z)$ as $\cyl :=
W_A \times (0 , z_{e})$, with $a = a_e$ and $\alpha = \theta_e$ as
above.

Let us consider first a smooth function $v$ on $W_e$ such that
$v(r, 0) = v(r, a) = 0$. We then have the one-dimensional
Poincar\'{e} inequality
\begin{equation}\label{eq.angle}
    \int_{0}^{\theta_{e}} |v|^2\; d\theta \le
    \frac{\pi}{\theta_{e}} \int_{0}^{\theta_{e}}
    |\partial_{\theta} v|^2\; d\theta.
\end{equation}
By integrating in polar coordinates we obtain
\begin{equation}\label{eq.poincare.cyl}
    \int_{W_a} \frac{|u|^2}{r^2}\; dx dy =
    \int_{W_a} \frac{|u|^2}{r}\; dr d\theta
    \le \frac{\pi}{\theta_{e}} \int_{W_a}
    \left ( \frac{|\partial_\theta u|^2}{r} \right )
    dr d\theta.
\end{equation}

We now claim that any $u \in C^{\infty}(\cyl)$ satisfying $u(r, 0,
z) = 0$ is such that
\begin{equation}\label{eq.I}
    \int_{\cyl} \frac{|u(x)|^2}{r^2} dV \le
    \frac{\pi}{\theta_{e}} \int_{\cyl} |\nabla u(x)|^2 dV \le
    \frac{\pi}{\theta_{e}} \int_{\Omega} |\nabla u(x)|^2 dV.
\end{equation}
Indeed, using Equation \eqref{eq.poincare.cyl} and the formula for
the $|\nabla u|$ in cylindrical coordinates, we get
\begin{multline*}
    \int_{\cyl} \frac{|u|^2}{r^2} dV =
    \int_{0}^{z_e} \int_{W_a } \frac{|u|^2}{r}\; dr d\theta dz
    \le \frac{\pi}{\theta_{e}} \int_{0}^{z_e} \int_{W_a}
    \left ( \frac{|\partial_\theta u|^2}{r} \right )
    dr d\theta dz \\
    \le  \frac{\pi}{\theta_{e}} \int_{0}^{z_e} \int_{W_a}
    \left ( \frac{|\partial_\theta u|^2}{r} +
    r |\partial_\rho u|^2 + r |\partial_z u|^2 \right )
    dr d\theta dz = \frac{\pi}{\theta_{e}} \int_{\cyl} |\nabla u(x)|^2 dV.
\end{multline*}

If $\cyl = \Omega_{v, e}$, then the proof proceeds exactly in the
same way, except that we replace $W_a$ in the integrals of the
last equation with $W_{za}$.

Now, if $\cyl = \Omega_v$, we proceed as in Equation
\eqref{eq.angle}, using also the formula
\begin{equation}
    |\nabla u|^2 = u^2_\rho + \frac 1 {\rho^2} u^2_\phi + \frac 1
    {\rho^2 \sin^2 \phi } u^2_\theta
\end{equation}
and the relation
\begin{equation*}
    \int_{\omega_v} |u|^2 \,dS
    \le C_v \int_{\omega_v} \left ( u^2_\phi + \frac 1
    {\sin^2 \phi } u^2_\theta \right) \sin \phi \,d\phi d\theta
    = C_v \int_{\omega_v} |\nabla_{(\phi, \theta)} u|^2 \,dS,
\end{equation*}
which is nothing but the usual Poincar\'{e}'s inequality on
$\omega_v$ ($dS$ is the volume element on $\omega_v$). We then
obtain ($\cyl = \Omega_v$)
\begin{multline}\label{eq.III}
    \int_{\cyl} \frac{|u|^2}{\rho^2} dV =
    \int_{0}^{2\epsilon}\! \int_{\omega_v } |u|^2 dS
    d\rho \le C_v \int_{0}^{2 \epsilon} \int_{\omega_v}
    \left ( u^2_\phi + \frac 1 {\sin^2 \phi } u^2_\theta \right )
    dS d\rho \\
    \le C_v \int_{0}^{2 \epsilon} \!\int_{\omega_v}
    \left (  u^2_\rho + \frac 1 {\rho^2} u^2_\phi + \frac 1
    {\rho^2 \sin^2 \phi } u^2_\theta \right )
    \rho^2 dS d\rho  = C_v
    \int_{\cyl} |\nabla u(x)|^2 dV.
\end{multline}

Adding the inequalities \eqref{eq.I} for all $\cyl = \Omega_{e}$
and all $\cyl = \Omega_{v, e}$, the inequalities \eqref{eq.III}
for all $\cyl = \Omega_{v}$, and the usual Poincar\'{e}
inequality, $\int_{\Omega} |u|^2 dV \le C\int_{\Omega} |\nabla
u(x)|^2 dV$, we obtain
\begin{equation*}
    \int_{\Omega} h |u|^2 dV \le
    C \int_{\Omega} |\nabla u(x)|^2  dV,
\end{equation*}
where $h(x)$ is a sum of $1$ and terms of the form $r^{-2}$,
$\rho^{-2}$. There will be one term $r^{-2}$ each time when $x \in
\Omega_{e}$ or $x \in \Omega_{v, e}$ and one term $\rho^{-2}$ each
time when $x \in \Omega_{v}$. Therefore $h \ge C r^{-2}$ on
$\Omega_{e}$ and on $\Omega_{v, e}$, $h \ge C \rho^{-2} \ge
Cr^{-2}$ on $\Omega_{v}$ and {\em outside} all of $\Omega_{e} \cup
\Omega_{v, e}$, and, finally, $h \ge 1 \ge C r^{-2}$ outside of
$\Omega_{e} \cup \Omega_{v, e} \cup \Omega_v$. Therefore $h \ge C
r^{-2}$ on {\em the whole of} $\Omega$. This completes the proof
of our inequality for $u$ smooth.
By a standard density argument, we obtain the desired result
\eqref{eq.newpoincare} for all functions in $H^1_0(\Omega)$.

\section{Proof of the main result\label{sec.Main}}

In this section, we prove our main result, Theorem
\ref{theorem.Main0}, assuming the intermediate results stated in the
previous sections.  We shall follow the pattern of proof from
\cite{Evans}. First, let us notice that Theorem \ref{thm.Sobolev}
allows us to reduce the proof to the case when $g = 0$.

Recall the function $r_{\Omega}$ introduced at the end of Section
\ref{sec.PM}. We can check directly that $r_{\Omega}^{\lambda}
\Delta r_{\Omega}^{-\lambda}$ depends continuously on $\lambda$
and that $r_{\Omega}^t\Kond{\mu}{a}(\Omega) =
\Kond{\mu}{a+t}(\Omega)$ (see \cite{BNZ2} for details in the case
of higher dimensions). This allows us to reduce the proof to the
case $a = 0$, as in \cite{BNZ1}.

We shall denote by $\big( u, v \big) := \int_{\Omega} u(x) v(x)
dx$ the (real) inner product on $L^2(\Omega)$. Let $\maH \subset
\Kond{1}{1}(\Omega)$ be the subspace consisting of the functions
$u \in \Kond{1}{1}(\Omega)$ such that $u = 0$ on $\pa \Omega$.
Thus $\maH$ is the kernel of the trace map $\Kond{1}{1}(\Omega)
\to \Kond{1/2}{1/2}(\pa \Omega)$. The Hardy--Poincar\'{e}
inequality (Theorem \ref{thm.Poincare}) then gives the following
inequality
\begin{equation*}
    B(u, u) := \big( \Delta u, u \big) = \sum_{j = 1}^n
    \big(\pa_j u, \pa_j u \big) \ge \epsilon
    \|u\|_{\Kond{1}{1}(\Omega)}^2,
\end{equation*}
for any $u \in \maH$. In particular, $B$ defines a continuous,
bilinear, coercive form on $\maH$. The assumptions of the
Lax-Milgram lemma are therefore satisfied, and hence $\Delta :
\maH \to \maH^* = \Kond{-1}{-1}(\Omega)$ is an isomorphism (by
this we understand that
$\Delta$ is continuous with continuous inverse), by the definition
of negative order Sobolev spaces on $\Omega$. This proves the
result for $\mu = 0$.

We now prove the result for an arbitrary $\mu \in \ZZ_+$.
Theorem \ref{thm.reg} and the result we
have just proved for $\mu = 0$ give that the map
\begin{equation}
    \Delta : \Kond{\mu+2}{1}(\Omega) \cap \{u \vert_{\pa \Omega}
    = 0\} \to \Kond{\mu}{-1}(\Omega)
\end{equation}
is surjective. Since this map is also continuous (Proposition
\ref{prop.map}) and injective (from the case $\mu = 0$), it is an
isomorphism by the open mapping theorem. The proof is now
complete.


\def\cprime{$'$} \def\cprime{$'$} \def\cprime{$'$} \def\cprime{$'$}
  \def\cprime{$'$} \def\cprime{$'$} \def\cprime{$'$} \def\cprime{$'$}
  \def\ocirc#1{\ifmmode\setbox0=\hbox{$#1$}\dimen0=\ht0 \advance\dimen0
  by1pt\rlap{\hbox to\wd0{\hss\raise\dimen0
  \hbox{\hskip.2em$\scriptscriptstyle\circ$}\hss}}#1\else {\accent"17 #1}\fi}

\end{document}